\documentclass[12pt,reqno]{article}

\usepackage[usenames]{color}
\usepackage{amssymb}
\usepackage{amsmath}
\usepackage{amsthm}
\usepackage{amsfonts}
\usepackage{amscd}
\usepackage{graphicx}

\usepackage{tikz}
\usepackage[colorlinks=true,
linkcolor=webgreen,
filecolor=webbrown,
citecolor=webgreen]{hyperref}

\definecolor{webgreen}{rgb}{0,.5,0}
\definecolor{webbrown}{rgb}{.6,0,0}

\usepackage{color}
\usepackage{fullpage}
\usepackage{float}

\usepackage{graphics}
\usepackage{latexsym}
\usepackage{epsf}
\usepackage{breakurl}

\begin{document}
	
\begin{center}
\epsfxsize=4in
\end{center}

	\theoremstyle{plain}
	\newtheorem{theorem}{Theorem}
	\newtheorem{corollary}[theorem]{Corollary}
	\newtheorem{lemma}[theorem]{Lemma}
	\newtheorem{proposition}[theorem]{Proposition}
	
	\theoremstyle{definition}
	\newtheorem{definition}[theorem]{Definition}
	\newtheorem{example}[theorem]{Example}
	\newtheorem{conjecture}[theorem]{Conjecture}
	
	\theoremstyle{remark}
	\newtheorem{remark}[theorem]{Remark}
	
	\begin{center}
		\vskip 1cm{\LARGE\bf 
			Skew Dyck paths having no peaks at level 1 (Sequence A128723)
		}
		\vskip 1cm
		\large
Helmut Prodinger \\
Stellenbosch University \\
		Department of Mathematical Sciences\\
		7602 Stellenbosch \\
		South Africa, \\
		and \\ NITheCS\\(National Institute for Theoretical and Computational Sciences)\\
		South Africa\\
		{\tt hproding@sun.ac.za} \\
	\end{center}
	
	\vskip .2 in
	
	\begin{abstract} 
		Skew Dyck paths are a variation of Dyck paths, where additionally to steps $(1,1)$ and $(1,-1)$ a south-west step $(-1,-1)$ is also allowed, provided that the path
		does not intersect itself. Replacing the south-west step by a red south-east step, we end up with decorated Dyck paths.
		Sequence A128723 of the Encyclopedia of Integer Sequences considers such paths where peaks at level 1 are forbidden. We provide a thorough analysis of a more general scenario, namely partial
		decorated Dyck paths, ending on a prescribed level $j$, both from left-to-right and from right-to-left (decorated Dyck paths are not symmetric).
		The approach is completely based on generating functions.
	\end{abstract}

\begin{center}
	\end{center}
	
	\section{Introduction}
	
			Skew Dyck paths are a variation of Dyck paths, where additionally to steps $(1,1)$ and $(1,-1)$ a south-west step $(-1,-1)$ is also allowed, provided that the path
	does not intersect itself. Replacing the south-west step by a red south-east step, we end up with decorated Dyck paths. Our earlier publication \cite{skew-old} studied such paths using generating functions.
	
	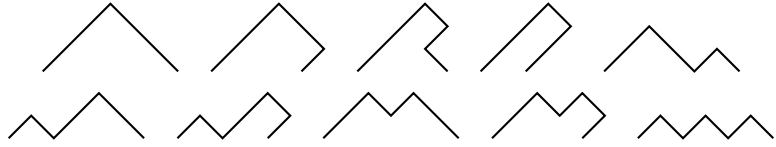
\begin{figure}[h]
		\begin{equation*}
			\begin{tikzpicture}[scale=0.3]
				\draw [thick](0,0)--(3,3)--(6,0);
			\end{tikzpicture}
			\quad
			\begin{tikzpicture}[scale=0.3]
				\draw [thick](0,0)--(3,3)--(5,1)--(4,0);
			\end{tikzpicture}
			\quad
			\begin{tikzpicture}[scale=0.3]
				\draw [thick](0,0)--(3,3)--(4,2)--(3,1)--(4,0);
			\end{tikzpicture}
			\quad
			\begin{tikzpicture}[scale=0.3]
				\draw [thick](0,0)--(3,3)--(4,2)--(3,1)--(2,0);
			\end{tikzpicture}
			\quad
			\begin{tikzpicture}[scale=0.3]
				\draw [thick](0,0)--(2,2)--(4,0)--(5,1)--(6,0);
			\end{tikzpicture}
		\end{equation*}
		\begin{equation*}
			\begin{tikzpicture}[scale=0.3]
				\draw [thick](0,0)--(1,1)--(2,0)--(4,2)--(6,0);
			\end{tikzpicture}
			\quad
			\begin{tikzpicture}[scale=0.3]
				\draw [thick](0,0)--(1,1)--(2,0)--(4,2)--(5,1)--(4,0);
			\end{tikzpicture}
			\quad
			\begin{tikzpicture}[scale=0.3]
				\draw [thick](0,0)--(1,1)--(2,2)--(3,1)--(4,2)--(6,0);
			\end{tikzpicture}
			\quad
			\begin{tikzpicture}[scale=0.3]
				\draw [thick](0,0)--(1,1)--(2,2)--(3,1)--(4,2)--(5,1)--(4,0);
			\end{tikzpicture}
			\quad
			\begin{tikzpicture}[scale=0.3]
				\draw [thick](0,0)--(1,1)--(2,0)--(3,1)--(4,0)--(5,1)--(6,0);
			\end{tikzpicture}
		\end{equation*}
		\caption{All 10 skew Dyck paths of length 6 (consisting of 6 steps).}
		\label{F1}
	\end{figure}
	
	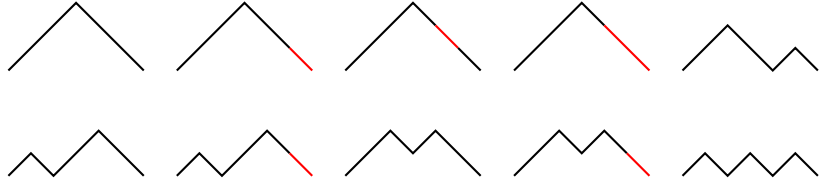
\begin{figure}[h]
		
		\begin{equation*}
			\begin{tikzpicture}[scale=0.3]
				\draw [thick](0,0)--(3,3)--(6,0);
			\end{tikzpicture}
			\quad
			\begin{tikzpicture}[scale=0.3]
				\draw [thick](0,0)--(3,3)--(5,1);
				\draw [thick,red](5,1)--(6,0);
			\end{tikzpicture}
			\quad
			\begin{tikzpicture}[scale=0.3]
				\draw [thick](0,0)--(3,3)--(4,2);
				\draw[red,thick] (4,2)--(5,1);
				\draw [thick](5,1)--(6,0);
			\end{tikzpicture}
\quad
			\begin{tikzpicture}[scale=0.3]
				\draw [thick](0,0)--(3,3)--(4,2);
				\draw[red,thick](4,2)--(6,0);
			\end{tikzpicture}
			\quad
			\begin{tikzpicture}[scale=0.3]
				\draw [thick](0,0)--(2,2)--(4,0)--(5,1)--(6,0);
			\end{tikzpicture}
		\end{equation*}
	
		\begin{equation*}
			\begin{tikzpicture}[scale=0.3]
				\draw [thick](0,0)--(1,1)--(2,0)--(4,2)--(6,0);
			\end{tikzpicture}
			\quad
			\begin{tikzpicture}[scale=0.3]
				\draw [thick](0,0)--(1,1)--(2,0)--(4,2)--(5,1);
				\draw[red,thick] (5,1)--(6,0);
			\end{tikzpicture}
			\quad
			\begin{tikzpicture}[scale=0.3]
				\draw [thick](0,0)--(1,1)--(2,2)--(3,1)--(4,2)--(6,0);
			\end{tikzpicture}
			\quad
			\begin{tikzpicture}[scale=0.3]
				\draw [thick](0,0)--(1,1)--(2,2)--(3,1)--(4,2)--(5,1);
				\draw[red,thick] (5,1)--(6,0);
			\end{tikzpicture}
			\quad
			\begin{tikzpicture}[scale=0.3]
				\draw [thick](0,0)--(1,1)--(2,0)--(3,1)--(4,0)--(5,1)--(6,0);
			\end{tikzpicture}
		\end{equation*}
		\caption{The 10 paths redrawn, with red south-east edges instead of south-west edges.}
		\label{F2}
	\end{figure}

Sequence A128723 considers such paths where peaks at level 1 are forbidden. These paths are the main objects of the present paper.	
	
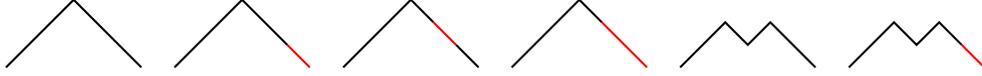
\begin{figure}[h]
		\begin{equation*}
			\begin{tikzpicture}[scale=0.3]
				\draw [thick](0,0)--(3,3)--(6,0);
			\end{tikzpicture}
			\quad
			\begin{tikzpicture}[scale=0.3]
				\draw [thick](0,0)--(3,3)--(5,1);
				\draw [thick,red](5,1)--(6,0);
			\end{tikzpicture}
			\quad
			\begin{tikzpicture}[scale=0.3]
				\draw [thick](0,0)--(3,3)--(4,2);
				\draw[red,thick] (4,2)--(5,1);
				\draw [thick](5,1)--(6,0);
			\end{tikzpicture}
			\quad
			\begin{tikzpicture}[scale=0.3]
				\draw [thick](0,0)--(3,3)--(4,2);
				\draw[red,thick](4,2)--(6,0);
			\end{tikzpicture}
			\quad
			\begin{tikzpicture}[scale=0.3]
				\draw [thick](0,0)--(1,1)--(2,2)--(3,1)--(4,2)--(6,0);
			\end{tikzpicture}
			\quad
			\begin{tikzpicture}[scale=0.3]
				\draw [thick](0,0)--(1,1)--(2,2)--(3,1)--(4,2)--(5,1);
				\draw[red,thick] (5,1)--(6,0);
			\end{tikzpicture}
		 		\end{equation*}
		\caption{The 6 paths without peaks on level 1.}
		\label{F3}
	\end{figure}

	The Figures \ref{F1}, \ref{F2}, \ref{F3} describe such paths of length 6.
	
	We catch the essence of a decorated Dyck path using a state-diagram (Fig.~\ref{arse4}):
	
	\begin{figure}[h]

		\begin{center}
			\begin{tikzpicture}[scale=1.5]
				\draw (0,0) circle (0.1cm);
				\fill (0,0) circle (0.1cm);
				
				\foreach \x in {0,1,2,3,4,5,6,7,8}
				{
					\draw (\x,0) circle (0.05cm);
					\fill (\x,0) circle (0.05cm);
				}
				
				\foreach \x in {0,1,2,3,4,5,6,7,8}
				{
					\draw (\x,-1) circle (0.05cm);
					\fill (\x,-1) circle (0.05cm);
				}
				
				\foreach \x in {0,1,2,3,4,5,6,7,8}
				{
					\draw (\x,-2) circle (0.05cm);
					\fill (\x,-2) circle (0.05cm);
				}
				
				\foreach \x in {0,1,2,3,4,5,6,7}
				{
					\draw[ thick,-latex] (\x,0) -- (\x+1,0);
					
				}

				\foreach \x in {1,2,3,4,5,6,7}
				{
					\draw[thick,  -latex] (\x+1,0) to[out=200,in=70]  (\x,-1);

				}
				\draw[ultra thick,purple,  dashed,  -latex] (1,0) to[out=200,in=70]  (0,-1);

				\foreach \x in {0,1,2,3,4,5,6,7}
				{
					
					\draw[thick,  -latex] (\x,-1) to[out=30,in=250]  (\x+1,0);	
					
				}

				\foreach \x in {0,1,2,3,4,5,6,7}
				{
					\draw[ thick,-latex] (\x+1,-1) -- (\x,-1);
					
				}
				\foreach \x in {0,1,2,3,4,5,6,7}
				{
					\draw[ thick,-latex,red] (\x+1,-1) -- (\x,-2);
					
				}
				
				\foreach \x in {0,1,2,3,4,5,6,7}
				{
					\draw[ thick,-latex,red] (\x+1,-2) -- (\x,-2);
					
				}
				
				\foreach \x in {0,1,2,3,4,5,6,7}
				{
					\draw[ thick,-latex] (\x+1,-2) -- (\x,-1);
					
				}

			\end{tikzpicture}
		\end{center}
		\caption{Three layers of states according to the type of steps leading to them (up, down-black, down-red).}
		\label{arse4}
	\end{figure}
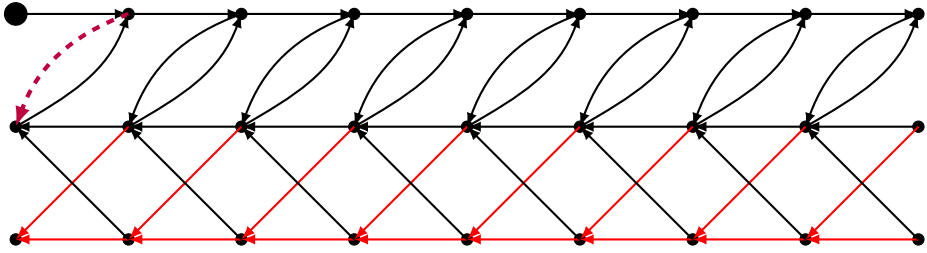

	It has three types of states, with $j$ ranging from 0 to infinity; in the drawing, only $j=0..8$ is shown. The first
	layer of states refers to an up-step leading to a state, the second layer refers to a black down-step leading to a state
	and the third layer refers to a red down-step leading to a state. 
	
	If the dashed edge is present, the graph models skew (decorated) Dyck paths. Any path from the origin to a node on level $j$
	represents such a decorated Dyck path ending on level $j$. In particular, if $j=0$, the path comes back to the $x$-axis.
	Note that the syntactic rules of forbidden patterns
	\begin{tikzpicture}[scale=0.3]\draw [thick](0,0)--(1,1); \draw [red,thick] (1,1)--(2,0);\end{tikzpicture}
	and
	\begin{tikzpicture}[scale=0.3] \draw [red,thick] (0,1)--(1,0);\draw [thick](1,0)--(2,1);\end{tikzpicture}
	can be clearly seen from the picture.
	
	However, if the dashed edge is \emph{not} present, it means that peaks at level 1 cannot be modeled by this graph, and that
	is what we want in the present paper.
	
	We will work out generating functions describing all paths
	leading to a particular state. We will use the notations $f_j,g_j,h_j$ for the three respective layers, from top to bottom.
	Although one could in principle compute all these functions separately, we are mainly interested in $s_j=f_j+g_j+h_j$, so
	we are interested in paths arriving on level $j$ but we do not care in which way this final level has been reached.
	It is also clear that a path of length $n$ leading to a state at level $j$ must satisfy $n\equiv j\bmod2$.
	
	In a last section, the right-to-left model is briefly described. Then, red down-steps become blue up-steps.

	\section{Generating functions and the kernel method}
	
	The functions depend on the variable $z$ (marking the number of steps), but mostly
	we just write $f_j$ instead of $f_j(z)$, etc.
	
	The following recursions can be read off immediately from the diagram (Fig. \ref{arse4}):
	\begin{gather*}
		f_0=1,\quad f_{i+1}=zf_i+zg_i,\quad i\ge0,\\
		g_i=zf_{i+1}+zg_{i+1}+zh_{i+1},\quad i\ge1,\\
		g_0=zg_{i+1}+zh_{i+1},\quad i\ge1,\\
				h_i=zh_{i+1}+zg_{i+1},\quad i\ge0.
	\end{gather*}
We can make a few direct observations: 	$f_0=1$, $f_1=z+zg_0$, $g_0=h_0$. The latter can be proved from combinatorial reasoning, since switching the last step from black to red resp.\ from red to black constitutes a bijection. This is a consequence of the fact that there are no peaks at level 1, otherwise the syntactic restrictions might be violated.

Now it is time to introduce  \emph{bivariate} generating functions: 
	\begin{equation*}
		F(z,u)=\sum_{i\ge0}f_i(z)u^i,\quad
		G(z,u)=\sum_{i\ge0}g_i(z)u^i,\quad
		H(z,u)=\sum_{i\ge0}h_i(z)u^i.
	\end{equation*}
	Again, often we just write $F(u)$ instead of $F(z,u)$ and treat $z$ as a `silent' variable. Summing the recursions leads to
	\begin{align*}
		\sum_{i\ge0}u^if_{i+1}&=\sum_{i\ge0}u^izf_i+\sum_{i\ge0}u^izg_i,\\
		\sum_{i\ge1}u^ig_i&=\sum_{i\ge1}u^izf_{i+1}+\sum_{i\ge1}u^izg_{i+1}+\sum_{i\ge1}u^izh_{i+1},\\
		\sum_{i\ge0}u^ih_i&=\sum_{i\ge0}u^izh_{i+1}+\sum_{i\ge0}u^izg_{i+1}.
	\end{align*}
	This can be rewritten as
	\begin{align*}
		\frac1u(F(u)-1)&=zF(u)+zG(u),\\*
		\sum_{i\ge1}u^ig_i+g_0&=\sum_{i\ge1}u^izf_{i+1}+\sum_{i\ge1}u^izg_{i+1}+zg_{1}+\sum_{i\ge1}u^izh_{i+1}+zh_1,\\
		\sum_{i\ge0}u^ig_i&=\sum_{i\ge1}u^izf_{i+1}+\sum_{i\ge0}u^izg_{i+1}+\sum_{i\ge0}u^izh_{i+1},\\
		G(u)&=\frac zu[F(u)-f_0-uf_1]+\frac zu[G(u)-g_0]+\frac zu[H(u)-h_0],\\
		H(u)&=	\frac zu(G(u)-G(0))+\frac zu(H(u)-H(0)).
	\end{align*}
Instead of working with 3 functions, we can reduce the system to just one equation (with the variable $G$):	
	\begin{equation*}
		F=\frac{1+zuG}{1-zu},\quad H=\frac{z(G-g_0-h_0)}{u-z}.
	\end{equation*}
Using this, we get
	\begin{equation*}
		G=\frac{-z^3u(u-z)+z(1-zu)(2+zu-z^2)g_0}{z(u-r_1)(u-r_2)}
	\end{equation*}
with
	\begin{equation*}
	r_1=\frac{1+z^2+\sqrt{1-6z^2+5z^4}}{2z},\quad r_2=\frac{1+z^2-\sqrt{1-6z^2+5z^4}}{2z}.
\end{equation*}
Note that $r_1r_2=2-z^2$. 
Since the factor $u-r_2$ in the denominator is ``bad,'' it must also cancel in the numerator. This is an essential step in
the kernel method, see for instance our own survey \cite{Prodinger-kernel}. This leads to the new version
	\begin{equation*}
		G=\frac{-z^3(u-z+r_2)-z^2(-z^2+zu+1+zr_2)g_0}{z(u-r_1)}.
	\end{equation*}
	Plugging in $u=0$ and solving the equation
	\begin{equation*}
		G(z,0)=g_0=\frac{-z^3(-z+r_2)-z^2(-z^2+1+zr_2)g_0}{z(-r_1)}
	\end{equation*}
leads to 
	\begin{equation*}
		g_0={\frac {1-2{z}^{4}-3{z}^{2}-\sqrt { 1-6z^2+5z^4 }}{ 2( {z}^{2}+		3 ) {z}^{2}}}.
	\end{equation*}
Knowing this, we know $G$, and thus $F$ and $H$. As the first goal, we still set $u=0$, thus considering paths coming back to the $x$-axis. Using Maple, 
	\begin{equation*}
	s_0:=	f_0+g_0+h_0={\frac {1-{z}^{4}-\sqrt { 1-6z^2+5z^4 }}{ ( {z}^{2}+		3 ) {z}^{2}}}.
	\end{equation*}

	\subsection*{The conjecture}
	We write $z^2=x$, since skew paths, as discussed here, must have an even number of steps.
	The function
	\begin{equation*}
		y(x)={\frac {1-{x}^{2}-\sqrt { 1-6x+5x^2 }}{ x( x+		3 ) }}
	\end{equation*}
	is the generating function of the sequence A128723:
	\begin{equation*}
		1, 0, 2, 6, 22, 84, 334, 1368, 5734, 24480, 106086, 465462, 2063658, 9231084, 41610162, \dots
	\end{equation*}
	\textsc{Gfun}, as described in \cite{gfun}, produces the algebraic equation that $y(x)$ satisfies:
	\begin{equation*}
		-(x-1)(x-2)+3x+2(1-x^2)y-x(3+x)y^2=0
	\end{equation*}
and	from this the differential equation
	\begin{equation*}
		- \left( 9 {x}^{2}+5 {x}^{3}+3-17 x \right) xy' + \left( 9 {x}^{2}+7 x-5 {x}^{3}-3 \right) y
		+3+9 {x}^{2}-5 {x}^{3}-7 x
	\end{equation*}
	and finally from the differential equation the recursion for the coefficients $s_n=[x^n]y(x)$:
	\begin{equation*}
		3(n+4)s_{n+3}-(17n+41)s_{n+2}+9ns_{n+1}+5(n+1)s_n=0.
	\end{equation*}
	An equivalent recursion was conjectured in the description of sequence A128723~\cite{OEIS}.
	
\subsection*{Partial paths}	

Another compution with Maple leads to
\begin{equation*}
S(z,u)=F(z,u)+G(z,u)+H(z,u)={\frac {-{z}^{4}-{z}^{4}g_0-{z}^{2}g_0+{z}^{2}-1}{z(u-r_1)}}.
\end{equation*}
	Further
	\begin{align*}
s_j:=[u^j]S(z,u)&={\frac {{z}^{4}+{z}^{4}g_0+{z}^{2}g_0-{z}^{2}+1}{zr_1(1-u/r_1)}}\\
&={\frac {{z}^{4}+{z}^{4}g_0+{z}^{2}g_0-{z}^{2}+1}{zr_1^{j+1}}}.
	\end{align*}
	One sees the parity: $j$ even/odd iff exponents are even/odd.
	If it is desired, $1/r_1$ may be expressed by $r_2$ (and a factor).
	
	\subsection*{Open-ended paths}
	
	We might allow \emph{any} level as end-level of the path. In terms of generating functions, this means to consider $S(z,1)$, viz.
	\begin{equation*}
S(z,1)=\frac{-2 {z}^{5}-3{z}^{4}+{z}^{3}-5 {z}^{2}-3 z+4-(z^2+3z+4)\sqrt{1-6z^2+5z^4}}{2z(3+z^2)(z^2+2z-1)}.
	\end{equation*}
The sequence of coefficients
\begin{equation*}
1, 1, 1, 2, 5, 8, 18, 31, 71, 126, 290, 527, 1218, 2253, 5223, 9796, 22763, 43170, 100502, 192347, \dots
\end{equation*}
is not in the encyclopedia \cite{OEIS}.

\section{Reading the decorated paths from right to left}

Since decorated paths are not symmetric, it makes sense to consider this scenario separately.
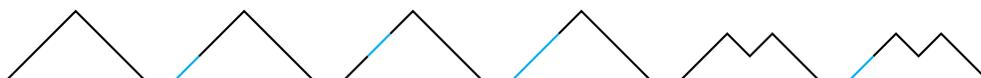
\begin{figure}[h]
	\begin{equation*}
		\begin{tikzpicture}[scale=0.3]
			\draw [thick](0,0)--(3,3)--(6,0);
		\end{tikzpicture}
		\quad
		\begin{tikzpicture}[scale=0.3]
			\draw [thick,cyan](0,0)--(1,1);
			\draw [thick](1,1)--(3,3)--(6,0);
		\end{tikzpicture}
		\quad
		\begin{tikzpicture}[scale=0.3]
			\draw [thick](0,0)--(1,1);
			\draw [thick,cyan](1,1)--(2,2);
			\draw [thick](2,2)--(3,3)--(6,0);
		\end{tikzpicture}
		\quad
		\begin{tikzpicture}[scale=0.3]
			\draw [thick,cyan](0,0)--(2,2);
			\draw [thick](2,2)--(3,3)--(6,0);
		\end{tikzpicture}
		\quad
						\begin{tikzpicture}[scale=0.3]
			\draw [thick](0,0)--(1,1)--(2,2)--(3,1)--(4,2)--(6,0);
		\end{tikzpicture}
		\quad
		\begin{tikzpicture}[scale=0.3]
			\draw [thick,cyan](0,0)--(1,1);
			\draw [thick](1,1)--(2,2)--(3,1)--(4,2)--(6,0);
		\end{tikzpicture}
			\end{equation*}
	\caption{All 6 dual (=right-to-left) skew Dyck paths of length 6 (consisting of 6 steps), having no peak at level 1.}
\end{figure}

We catch the essence of a decorated (dual skew) Dyck path using a state-diagram:

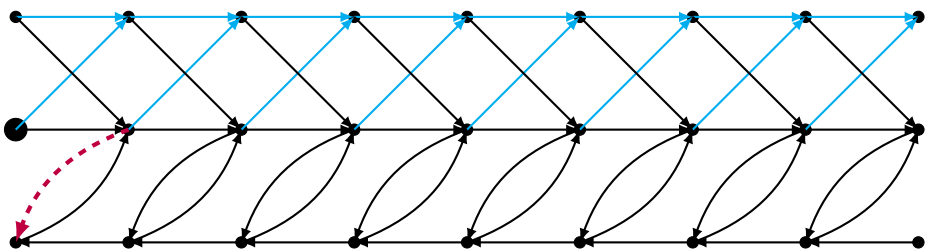
\begin{figure}[h]
	
	\begin{center}
		\begin{tikzpicture}[scale=1.5]
			\draw (0,0) circle (0.1cm);
			\fill (0,0) circle (0.1cm);
			
			\foreach \x in {0,1,2,3,4,5,6,7,8}
			{
				\draw (\x,0) circle (0.05cm);
				\fill (\x,0) circle (0.05cm);
			}
			
			\foreach \x in {0,1,2,3,4,5,6,7,8}
			{
				\draw (\x,-1) circle (0.05cm);
				\fill (\x,-1) circle (0.05cm);
			}
			
			\foreach \x in {0,1,2,3,4,5,6,7,8}
			{
				\draw (\x,1) circle (0.05cm);
				\fill (\x,1) circle (0.05cm);
			}
			
			\foreach \x in {0,1,2,3,4,5,6,7}
			{
				\draw[ thick,-latex] (\x,0) -- (\x+1,0);
				
			}
			\foreach \x in {0,1,2,3,4,5,6,7}
			{
				\draw[ thick,-latex,cyan] (\x,1) -- (\x+1,1);
				
			}
			
			\foreach \x in {0,1,2,3,4,5,6,7}
			{
				\draw[ thick,-latex,cyan] (\x,0) to (\x+1,1);
				
			}
			
			\foreach \x in {0,1,2,3,4,5,6,7}
			{
				\draw[ thick,-latex] (\x,1)  to (\x+1,0);
				
			}
			
			\foreach \x in {0,1,2,3,4,5,6,7}
			{
				\draw[ thick,latex-] (\x+1,0)   to[out=250,in=20](\x,-1);
				
			}
			
			\foreach \x in {0,1,2,3,4,5,6,7}
			{
				\draw[ thick,-latex] (\x+1,-1)  to (\x,-1);
				
			}
				\foreach \x in {1,2,3,4,5,6,7}
{
	\draw[thick,  -latex] (\x+1,0) to[out=200,in=70]  (\x,-1);

}
\draw[ultra thick,purple,  dashed,  -latex] (1,0) to[out=200,in=70]  (0,-1);

		\end{tikzpicture}
	\end{center}
	\caption{Three layers of states according to the type of steps leading to them (down, up-black, up-blue).}
\end{figure}
Note that the syntactic rules of forbidden patterns
\begin{tikzpicture}[scale=0.3]\draw [thick,cyan](0,0)--(1,1); \draw [thick] (1,1)--(2,0);\end{tikzpicture}
and
\begin{tikzpicture}[scale=0.3] \draw [thick] (0,1)--(1,0);\draw [thick,cyan](1,0)--(2,1);\end{tikzpicture}
can be clearly seen from the picture.

As in the earlier section, if the dashed edge is removed it means that the condition `no peak at level 1' is modeled,
which is what we need to do.
Using the letters $c_j,a_j,b_j$ (in that order) for the generating functions to reach state $j$ in the particular layer, we find the 
 following recursions  immediately from the diagram:
\begin{gather*}
	a_0=1,\quad a_{i+1}=za_i+zb_i+zc_i,\quad i\ge0,\\
	b_0=zb_1,\quad	b_i=za_{i+1}+zb_{i+1},\quad i\ge1,\\
	c_{i+1}=za_{i}+zc_{i},\quad i\ge0.
\end{gather*}
We introduce  \emph{bivariate} generating functions: 
\begin{equation*}
	A(z,u)=\sum_{i\ge0}a_i(z)u^i,\quad
	B(z,u)=\sum_{i\ge0}b_i(z)u^i,\quad
	C(z,u)=\sum_{i\ge0}c_i(z)u^i.
\end{equation*}
 Summing the recursions leads to
\begin{align*}
	A(u)&=\sum_{i\ge0}u^ia_i =1+u\sum_{i\ge0}u^i(za_i+zb_i+zc_i)\\
	&=1+uzA(u)+uzB(u)+uzC(u),\\
	\sum_{i\ge0}u^ib_i &= \sum_{i\ge1}u^iza_{i+1}+\sum_{i\ge0}u^izb_{i+1},\\
	B(u)&=\frac zu\sum_{i\ge2}u^ia_i+\frac zu\sum_{i\ge1}u^ib_i\\
	&=\frac zu[A(u)-a_0-ua_1]+\frac zu[B(u)-b_0],\\
	\sum_{i\ge1}u^ic_i &=uz\sum_{i\ge0}u^ia_i+uz\sum_{i\ge0}u^ic_i,\\
	C(u)-c_0&=uzA(u)+uzC(u).	
\end{align*}
We have $c_0=0$, $a_0=1$, and $a_1=z+zb_0$. We may write 
\begin{equation*}
C(u)=\frac{uzA(u)}{1-uz},
\end{equation*}
\begin{equation*}
	A(u)
=1+uzA(u)+uzB(u)+\frac{u^2z^2A(u)}{1-uz}=\frac{1+uzB(u)}{1-uz-\frac{u^2z^2}{1-uz}}=\frac{1-uz}{1-2uz}[1+uzB(u)],
\end{equation*}
\begin{equation*}
C(u)=\frac{uz}{1-2uz}[1+uzB(u)].
\end{equation*}
Solving for $B(u)$,
\begin{equation*}
B(u)=\frac{z \left( 2{u}^{2}{z}^{2}+2{z}^{2}{u}^{2}b_0+b_0zu-b_0 \right)}{z(z^2-2)(u-r_1^{-1})(u-r_2^{-1})}.
\end{equation*}
We cancel the bad factor $(u-r_1^{-1})$ out of the numerator:
\begin{equation*}
B(u)=\frac {z ( 2r_1uz+2r_1uzb_0+b_0r_1+2z+2zb_0 ) r_2}{r_1 ( {z}^{2}-	2 ) ( ur_2-1 ) }
\end{equation*}
Plugging in $u=0$ results in the equation
\begin{equation*}
b_0=\frac {-z ( b_0r_1+2z+2zb_0 ) r_2}{r_1 ( {z}^{2}-	2 )  }
\end{equation*}
and thus
\begin{equation*}
b_0=\frac{1-z^4-\sqrt{1-6z^2+5z^4}}{z^2(3+z^2)}-1,
\end{equation*}
as expected, since $1+b_0$ is the generating function of all skew Dyck paths without peaks at level 1.

Expressions for $A(z,u)+B(z,u)+C(z,u)$ and $[u^j](A(z,u)+B(z,u)+C(z,u))$ could be explicitly written, which we leave to the reader.

The open paths in this model are enumerated via
\begin{equation*}
A(z,1)+B(z,1)+C(z,1),
\end{equation*}
which is a long expression, with coefficients
\begin{equation*}
1, 2, 4, 10, 24, 56, 134, 318, 764, 1824, 4390, 10520, 25346, 60878, 146768,\dots
\end{equation*}
which are again not in the  OEIS \cite{OEIS}. 

Explicit formul\ae{} for this model are a bit unpleasant, but easily regenerated using Maple.

\section{Conclusion}

In order to keep this paper short (and not boring) we refrained from working out many additional parameters.
That might be a good project for graduate students.

	\bigskip
	\hrule
	\bigskip
	
	\noindent 2010 {\it Mathematics Subject Classification}:
	Primary 05A15; Secondary 05A19.
	
	\noindent \emph{Keywords:  Skew Dyck path, peak, forbidden pattern, generating function, kernel-method.}
	
	\bigskip
	\hrule
	\bigskip
	
	\noindent (Concerned with sequence
	{A128723}.)

	\vspace*{+.1in}
	\noindent
	
	\bigskip
	\hrule
	\bigskip
	
	\noindent
	Return to
	\htmladdnormallink{Journal of Integer Sequences home page}{https://cs.uwaterloo.ca/journals/JIS/}.
	\vskip .1in


\begin{thebibliography}{99}
		
	
		
	
\bibitem{Prodinger-kernel}
Helmut Prodinger.
\newblock {\em The Kernel Method: A Collection of Examples}, 
\newblock Séminaire Lotharingien de Combinatoire,  B50f (2004), 19 pp.

	\bibitem{skew-old} H. Prodinger, Partial skew Dyck paths---a kernel method approach, {\it arXiv:2108.09785}  (2021). 


\bibitem{gfun}
B.  Salvy and P.  Zimmermann.
\newblock {\textsc{Gfun}}: A Maple package for the manipulation of generating and holonomic functions in one variable.
\newblock {\em ACM Trans. Math. Softw.}, 20(2): 163--177 1994.



\bibitem{OEIS}
Neil J.~A. Sloane and The OEIS~Foundation Inc.
\newblock The on-line encyclopedia of integer sequences, 2022.


 


		
	\end{thebibliography}
\end{document}